\documentclass{Authors/styles/svmult}
\usepackage[utf8]{inputenc}
\usepackage[english]{babel}
\usepackage{amsfonts}
\usepackage{amsmath}
\usepackage{graphicx}
\usepackage{array}

\date\today

\newcommand{\eps}{{\epsilon}}

\begin{document}
\title*{Penalty Methods for the Hyperbolic System Modelling the Wall-Plasma Interaction in a Tokamak}
\author{Philippe Angot, Thomas Auphan, Olivier Gu\`es}
\institute{Philippe Angot, Thomas Auphan and Olivier Gu\`es \at Universit\'e de Provence, Laboratoire d'Analyse Topologie et Probabilit\'es, Centre de Math\'ematiques et Informatique, 39 rue Joliot Curie, 13453 Marseille Cedex 13, France, \email{[angot,tauphan,gues]@cmi.univ-mrs.fr}
}
\maketitle

\renewcommand{\b}[1]{\mathbf{#1}}
\renewcommand{\u}[1]{\underline{#1}}

\abstract{
 The penalization method is used to take account of obstacles in a tokamak, such as the limiter.  
We study a non linear hyperbolic system modelling the plasma transport in the area close to the wall. 
 A penalization which cuts the transport term of the momentum is studied. We show numerically that this penalization creates a Dirac measure at the plasma-limiter interface which prevents us from defining the transport term in the usual sense. Hence, a new penalty method is proposed for this hyperbolic system and numerical tests reveal an optimal convergence rate without any spurious boundary layer.
}

\section{Introduction}

A tokamak is a machine to study plasmas and the fusion reaction. The plasma at high temperature ($10^{8} K$) is confined in a toro\"idal chamber thanks to a magnetic field. One of the main goals is to perform controlled fusion with enough efficiency to be a reliable source of energy.
But, since the magnetic confinement is not perfect, the plasma is in contact with the wall. In order to preserve the integrity of the wall and to limit the pollution of the plasma, it is crucial to control these interactions.

We study, using a fluid approximation of the plasma, a simplified system of equations governing the plasma transport in the scrape-off layer, parallel to the magnetic field lines.
In this paper, after a numerical study of the pena\-lization introduced by Isoardi \emph{et al.} \cite{Iso10}, we modify the boundary conditions to ensure the well-posedness of the hyperbolic system and we propose another penalty method which seems to be free of boundary layer.

\section{The model hyperbolic problem}\label{Presentation of the hyperbolic system}
In this paper, we consider a very simple model taking only into account the transport in the direction parallel 
to the magnetic field lines, (see for example  \cite{Iso10}, \cite{Tam07}).
It is a one dimensional $2\times 2$ hyperbolic system of conservation laws 
for the particle density $N$ and the particle flux $\Gamma$, which reads :
\begin{equation}\label{O1}
\left\{
 \begin{array}{l}
     \partial_t N + \partial_x \Gamma=S \\
     \partial_t \Gamma + \partial_x \left(\dfrac{\Gamma^2}{N} + N\right)=0 \\
     \textsf{Initial conditions : }N(0,.)=N_0 \textsf{ and } \Gamma(0,.)=\Gamma_0
 \end{array}
 \right. (t,x) \in \mathbb{R}^+_* \times ]-L,L[ 
\end{equation}
Here, the boundaries of the domain  $x=L$ and $x=-L$
correspond to the "limiters", which are material obstacles for the fluid (see Fig. 1). In the right-hand side,
$S$ is a source term. 

There is a difficulty with the choice of the boundary conditions for the system (\ref{O1}).
From  physical arguments, it follows that the domain (namely the scrape-off layer) is basically divided
into two regions  \cite{Tam07} :
\begin{itemize}
 \item One region far from the limiter, the pre-sheath, where the plasma is neutral and the Mach number $M=\Gamma / N$ of the plasma satisfies $|M|\leq 1$.
 \item One region next to the limiter (in a thin layer called the sheath area, whose typical thickness is of the order of $10^{-5} m$), where the electroneutrality hypothesis does not hold and we have $|M|>1$. More precisely $M>1$ close to $x=L$ and
 $M<-1$ close to the boundary $x=-L$.
\end{itemize}
It could seem natural to prescribe $M=1$ (resp. $M=-1$)
as a boundary condition at $x=L$ (resp. $x=-L$) for the system,
since the physical arguments imply that $M=\pm 1$ very close to the obstacle (Bohm criterion).  
These are exactly the boundary conditions which are chosen in  \cite{Iso10}.
However, in that case, as the eigenvalues of the Jacobian of the flux function are $M-1$ and $M+1$, it follows that
at the plasma limiter interface one eigenvalue is $0$ (the boundary is characteristic) and the other one is outgoing (it is also true at $x=-L$), and clearly the problem does not satisfy the usual sufficient conditions for well posedness, see \cite{gue90}, \cite{rau85}, \cite{ben07} : the number of boundary conditions ($=1$) is not equal
to the number of incoming eigenvalues ($=0$). 

In order to test our penalty approach with a well-defined hyperbolic boundary value problem,
in section \ref{Study of penalty methods}, we slightly modify the boundary conditions
of the paper \cite{Iso10}, and impose $M= 1 -\epsilon$ on $x=L$ and $M= -1 + \epsilon$ on $x=-L$
with a fixed $\epsilon >0$, which leads to a well-posed hyperbolic problem. In our numerical simulations we use $\epsilon=0.1$.

\begin{figure}
\begin{center}
\includegraphics[scale=0.35, trim = 20mm 91mm 20mm 118mm, clip=true]{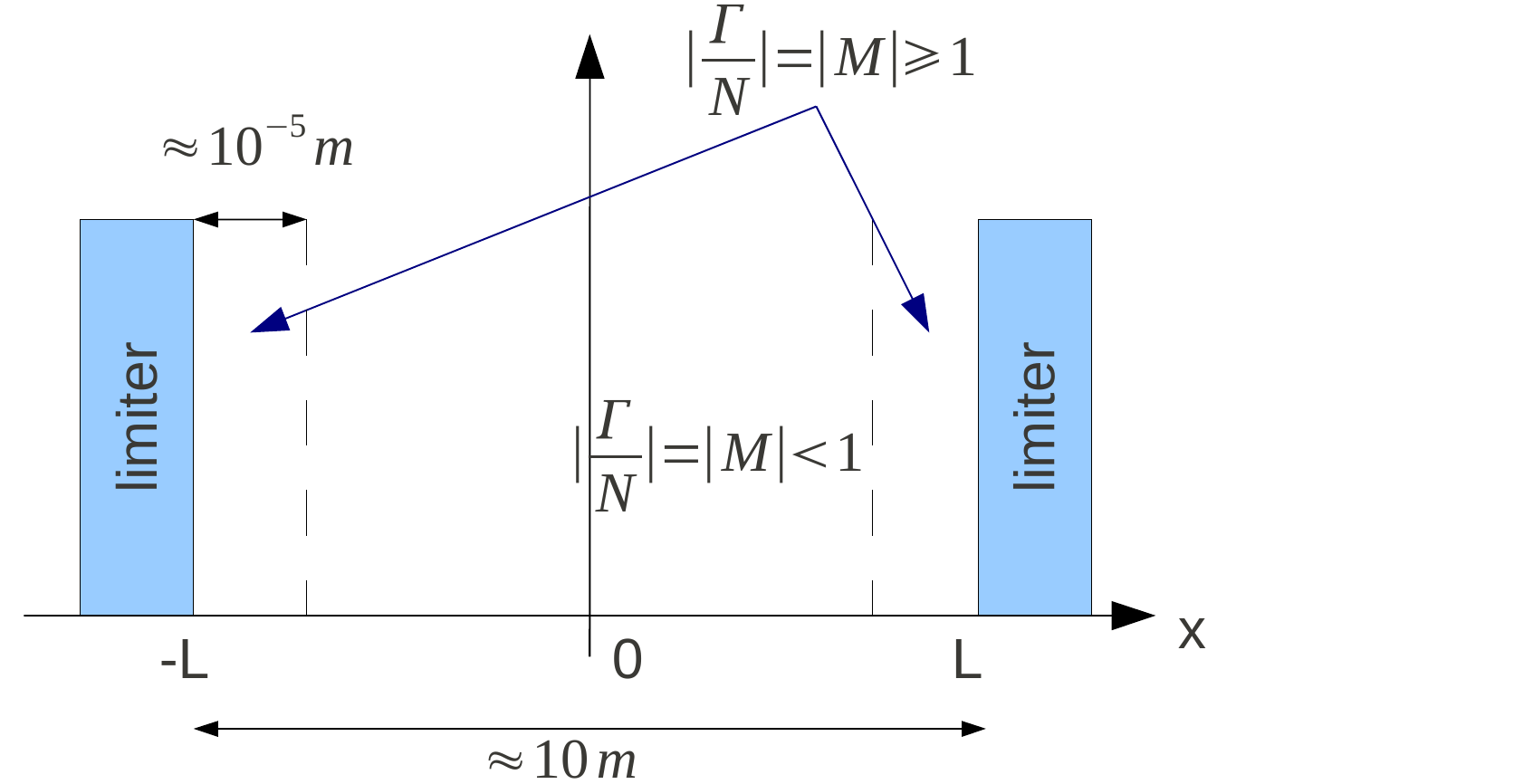}
\end{center}
\caption{
Schematic representation of the scrape-off layer. The $x$-axis corresponds to the curvilinear coordinate along a magnetic line close to the wall of the tokamak.
}
\label{Complete_domain}
\end{figure}

The numerical tests presented below, use a finite volume scheme with a second order extension : MUSCL reconstruction with the \emph{minmod} slope limiter and the Heun scheme which is a second order Runge--Kutta TVD time discretization. The finite volume scheme is the VFRoe using the non conservative variables for the linearized Riemann solver \cite{Gal03}; here, the non conservatives variables are $N$ and $M$. To avoid stability issues, the penalized terms are treated implicitly for the time discretization.

\section{Study of penalty methods}\label{Study of penalty methods}
\subsection{A first penalty method}

The following penalty approach has been proposed by Isoardi 
\emph{et al.} \cite{Iso10}. Let's $\chi$ be the characteristic function of the limiter, i.e. 
$\chi(x)=1$ if $x$ is in the limiter, and  $\chi(x)=0$ elsewhere, and
$\eta$ the penalization parameter. The penalized system is given by :
\begin{equation}\label{Penal_pb}
\left\{\begin{array}{l}
 \partial_t N + \partial_x \Gamma + \dfrac{\chi}{\eta} N = (1-\chi) S_N \qquad \textnormal{ in } \mathbb{R}^+_* \times \mathbb{R}\\
 \partial_t \Gamma + (1-\chi) \partial_x \left(\dfrac{\Gamma^2}{N} + N \right) +\dfrac{\chi}{\eta}(\Gamma-M_0 N)=(1-\chi) S_{\Gamma}\\
 \textsf{Initial conditions : }N(0,.)=N_0 \textsf{ and } \Gamma(0,.)=\Gamma_0
\end{array}\right.
\end{equation}
$M_0$ is a function such that, at the plasma-limiter interface we have $|M_0|=1$.
Here, the two components of the unknown are penalized although there is no incoming wave. At least formally, 
$N$ is forced to converge to $0$ inside the limiter when $\eta$ tends to $0$. 

The flux of the second equation is cut inside of the limiter, and this causes 
some troubles from the mathematical point of view. Indeed, this is
an hyperbolic system with discontinuous coefficients and the meaning of the term
$$
(1-\chi) \partial_x \left(\dfrac{\Gamma^2}{N} + N \right)
$$
is not clear because it can involve the product of a measure with a discontinuous function which has
no distributional sense. As a consequence and as a confirmation of this fact,
our numerical tests show the existence of a strong singularity at the interface for the numerical discrete solution. 
Concerning the interpretation of this numerical singularity, it could happen
(but we don't have any rigorous proof and this is just an open question) that this system admits 
generalized solutions in the spirit of Bouchut--James  \cite{Bou98} (see also Poupaud--Rascle \cite{Pou97}, or  Fornet--Gu\`es \cite{For08}) such as measure-valued solutions, which can for example exhibit a Dirac measure at the interface,
and this generalized solution could be selected by the numerical approximation process.  

For the numerical test, we choose $S_N$ and $S_{\Gamma}$ so that the following functions define a solution of the boundary value problem :
\begin{equation*}
 N(t,x)=\exp \left(\dfrac{-x^2}{0.16 (t+1)}\right) \qquad \Gamma(t,x)=\sin \left(\dfrac{\pi x}{0.8} \right) \exp \left(\dfrac{-x^2}{0.16 (t+1)}\right)
\end{equation*}
These test solutions are regular (at least inside the plasma area) and has no singularity at the plasma-limiter interface.
\begin{figure}
\begin{center}
\includegraphics[scale=0.56, trim = 10mm 8mm 10mm 11mm, clip=true]{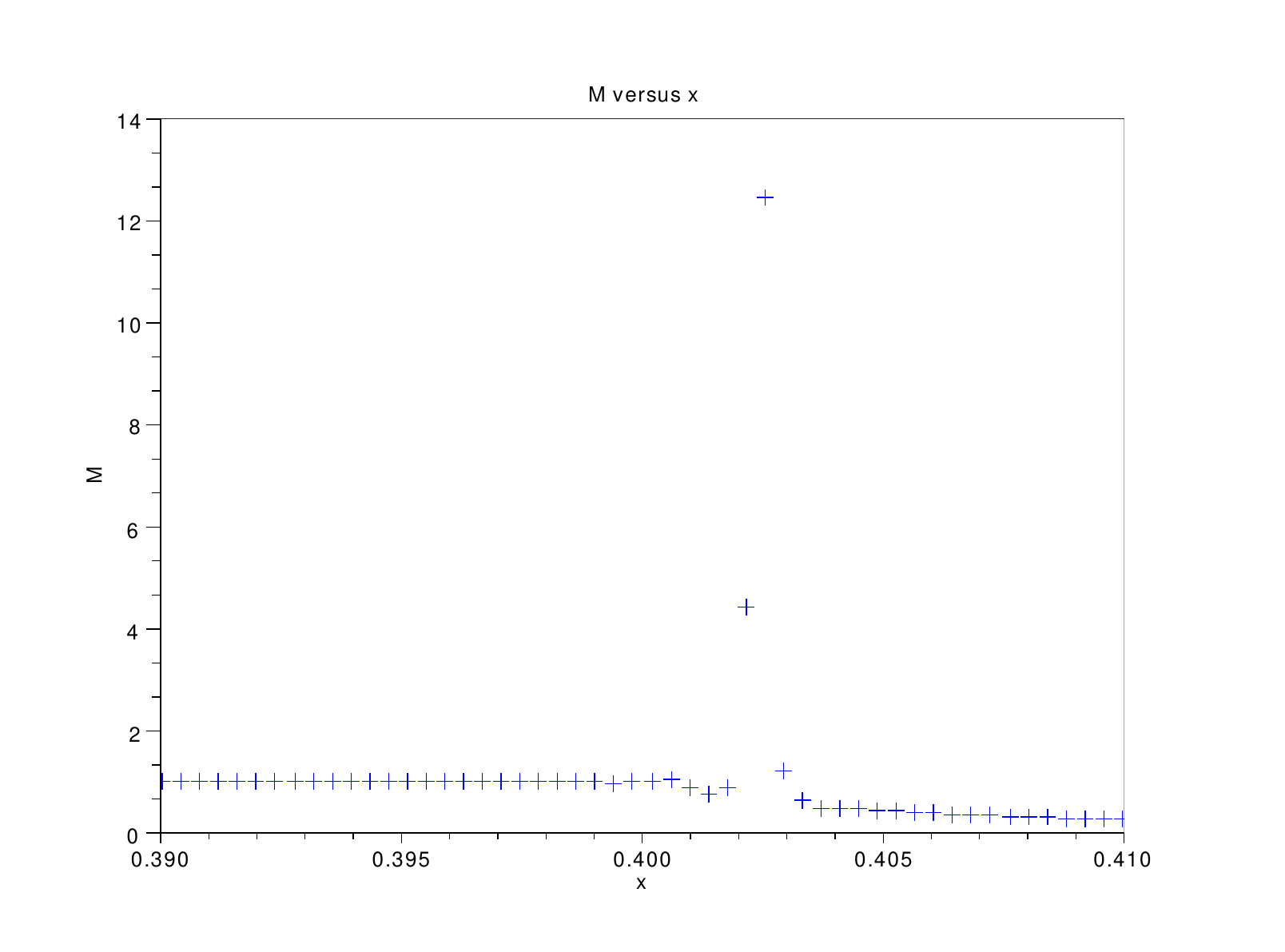}

\end{center}
\caption{
$M$ versus $x$ with $\eta=10^{-3}$, a mesh of $J = 1280$ cells using the penalization of Isoardi \emph{et al.} \cite{Iso10}. The computations are stopped when $\max_{i \in \{1,\dots,J\}}(|M_i^n|)>10$, which corresponds to the time : $t= 0.008822$. The computational domain was $[0,0.5]$ and $L=0.4$ (plasma-limiter interface). At $x=0$, we impose a symmetry condition.
}
\label{Dirac_measure}
\end{figure}
In the Fig. \ref{Dirac_measure}, we observe that a peak appears very quickly, then $|M_i^n|$ become very large (about $10^8$) in a few points.
The same computations are made for two more refined meshes (respectively for $2560$ and $10240$ cells) and we observe that the peak is nearer and nearer to the plasma limiter interface, when the resolution increases. Besides, when the mesh step decreases, the peak appears earlier and earlier. We stop the computations when $\max_{i \in \{1,\dots,J\}}(|M_i^n|)>10$ but similar results are obtained when the stop criterion is $\max_{i \in \{1,\dots,J\}}(|M_i^n|)>100$.
This leads one to believe that, if the solution converges to a generalized solution of the continuous problem, then 
this generalized solution must have a singularity supported by the interface (that could be a Dirac measure for example). 
We notice that the presence of a Dirac measure at the interface is not only a theoretical issue since it has been observed numerically and that the Dirac measure destabilizes numerical schemes. In the following section, we propose a modification of the boundary value problem to obtain a well-posed version.

\subsection{A new penalty method for the modified boundary conditions}

After the modifications proposed in section \ref{Presentation of the hyperbolic system}, the well-posed initial boundary value problem reads :
\begin{equation} \label{O2}
\left\{
 \begin{array}{l}
     \partial_t N + \partial_x \Gamma=S \\
       \partial_t \Gamma + \partial_x \left(\dfrac{\Gamma^2}{N} + N\right)=0 \\
       M(.,-L)=-1+\eps \textsf{ and } M(.,L)=1-\eps   \\
       N(0,.)=N_0 \textsf{ and } \Gamma(0,.)=\Gamma_0
    \end{array}
    \right .
    \qquad (t,x) \in \mathbb{R}^+_* \times ]-L,L[ 
\end{equation}
For this problem, the boundary is not characteristic, and the boundary conditions are
maximally dissipative. Hence, for compatible initial data, the problem has a unique local in time solution, which is regular enough : at least $\mathcal{C}^1$ is sufficient to perform the asymptotic analysis;
see e.g. \cite{ben07}, \cite{rau74}. 

To penalize (\ref{O2}), we use a method developed in the semi-linear case by Fornet and Gu\`es \cite{For09}. In order to have an homogeneous Dirichlet boundary condition for the theoretical study, the system is reformulated with the unknowns $\tilde{u}=\ln(N)$ and $\tilde{v}=\Gamma/N-M_0$. Although our system is quasi-linear (and not semi-linear), the method can be extended to this case. An interesting feature of the method is that it yields to a convergence
result without generation of a boundary layer inside the limiter. 
Up to now, we don't know if this method can be extended to more general quasi-linear first order hyperbolic system with maximally dissipative conditions.

We assume that $M_0$ is a constant such that $0<M_0<1$.
We denote by $\chi$ the characteristic function associated to the limiter, i.e. $\chi(x)=1$ if the point $x$ is in the limiter.

The new penalized problem reads :
\begin{equation}\label{Penal_pb_OK}
\left\{ \begin{array}{l}
 \partial_t N + \partial_x \Gamma = S_N\\
 \partial_t \Gamma + \partial_x \left(\dfrac{\Gamma^2}{N} + N\right) + \dfrac{\chi}{\eta}\left(\dfrac{\Gamma}{M_0} - N\right)=S_{\Gamma}\\
 N(0,.)=N_0 \textsf{ and } \Gamma(0,.)=\Gamma_0
 \end{array}\right. \qquad \textnormal{ in } \mathbb{R}^+_* \times \mathbb{R}
\end{equation}

The formal asymptotic expansion of a continuous solution to (\ref{Penal_pb_OK}) with the BKW (Brillouin--Kramers--Wentzel) method does not contain any boundary layer term \cite{Ang11} and this suggests strongly that there is no boundary layer at all in the solution.
Notice that the penalization is incomplete: only one field is penalized, which is natural since there is only one boundary condition.

For the numerical tests, we use a regular solution:
$$N(t,x)=\exp \left(\dfrac{-x^2}{0.16 (t+1)}\right) \qquad \Gamma(t,x)=M_0 \sin \left(\dfrac{\pi x}{0.8} \right) \exp \left(\dfrac{-x^2}{0.16 (t+1)}\right)$$
 and $S_N, S_{\Gamma}$ are well chosen. The spatial domain is $[0,0.5]$ with a symmetry condition at $x=0$ and the limiter set corresponds to $x \in [0.4,0.5]$.

\begin{figure}
\begin{center}
\includegraphics[scale=0.30, angle=-90, trim = 10mm 20mm 8mm 20mm, clip=true]{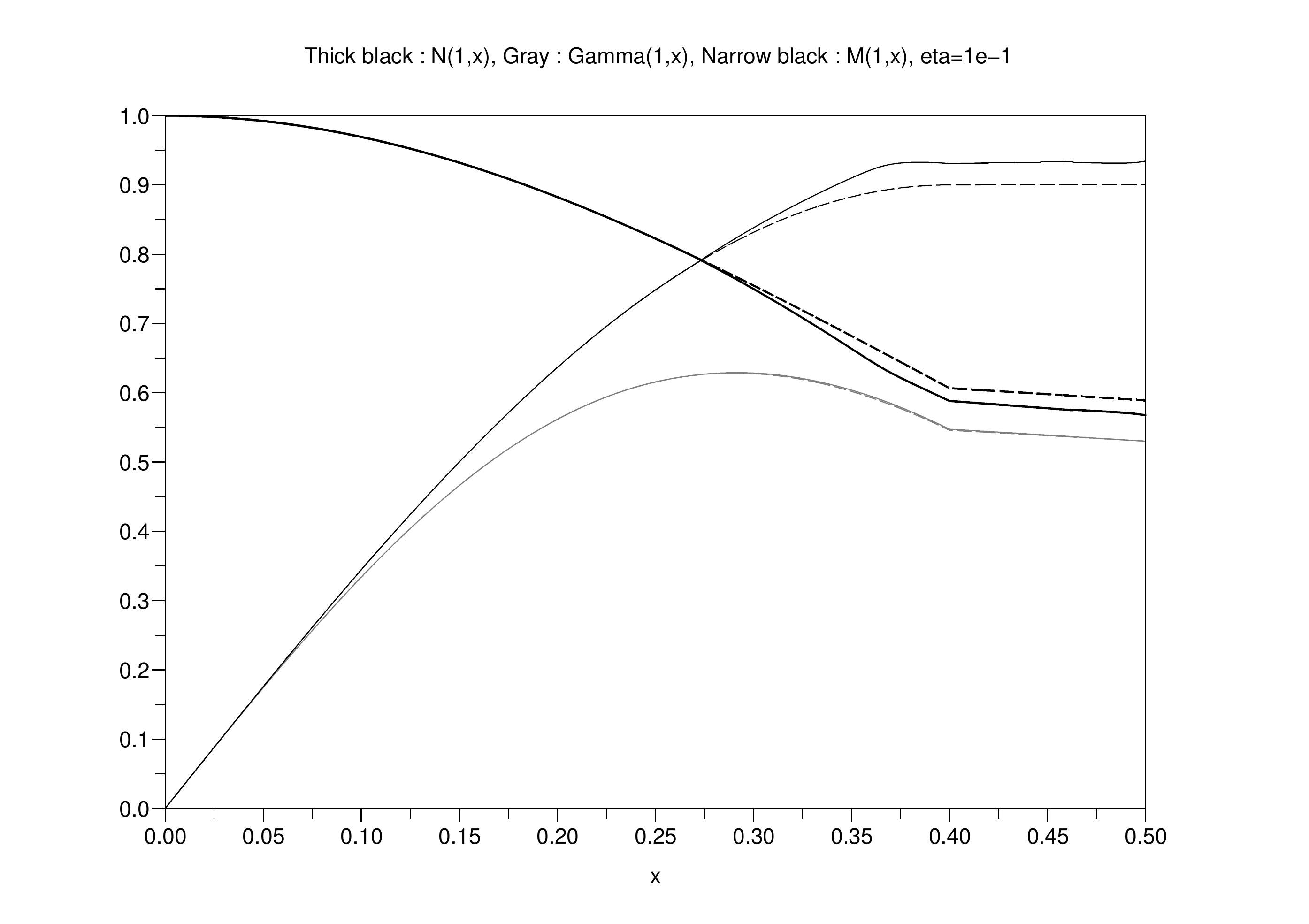}
\end{center}
\caption{
Plot of $N$, $\Gamma$ and $M$ as functions of $x$ (at $t=1$) with the penalty method free of boundary layer for $\eta=0.1$. The continuous lines represent the numerical solutions whereas the dashed lines corresponds to the exact solution of the hyperbolic limit problem ($\eta \to 0$). The limiter corresponds to the area $x\in[0.4,0.5]$. For smaller values of $\eta$, for instance for $\eta=10^{-5}$, the plot is almost the same as the plot of the exact solution (dotted lines).
}
\label{Plot_newpen_regular}
\end{figure}

We analyze the convergence when the penalization parameter $\eta$ tends to $0$ using a uniform spatial mesh of step $\delta x=10^{-5}$. We calculate the error in $L^1$ norm for $N$, $\partial_x N$, $\Gamma$ and $\partial_x \Gamma$. The goal is to confirm numerically the absence of boundary layer with an optimal rate of convergence as $\mathcal{O}(\eta)$.

\begin{figure}
\begin{center}
\includegraphics[scale=0.56, trim = 5mm 3mm 2mm 10mm, clip=true]{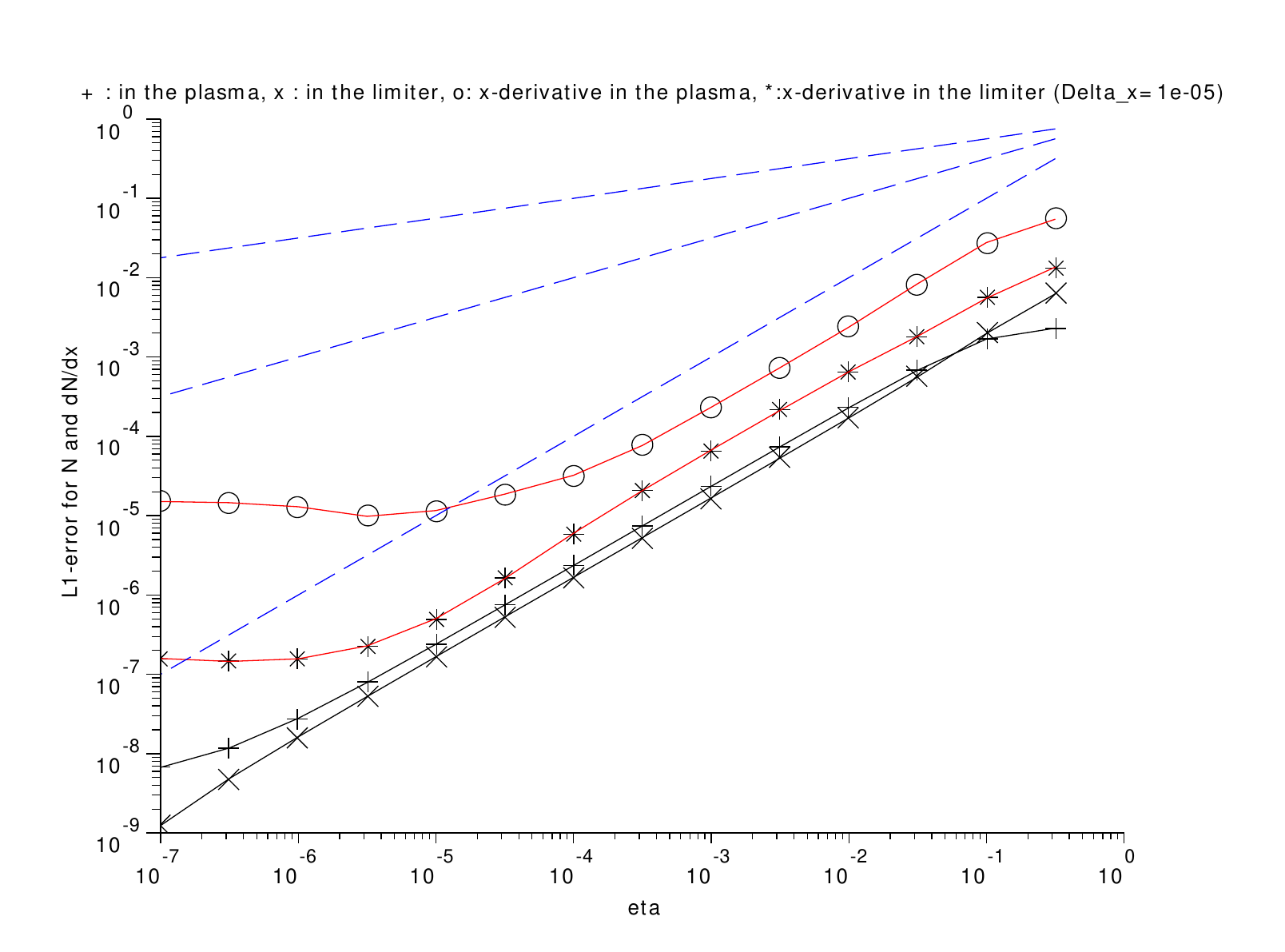}
\includegraphics[scale=0.56, trim = 5mm 4mm 2mm 10mm, clip=true]{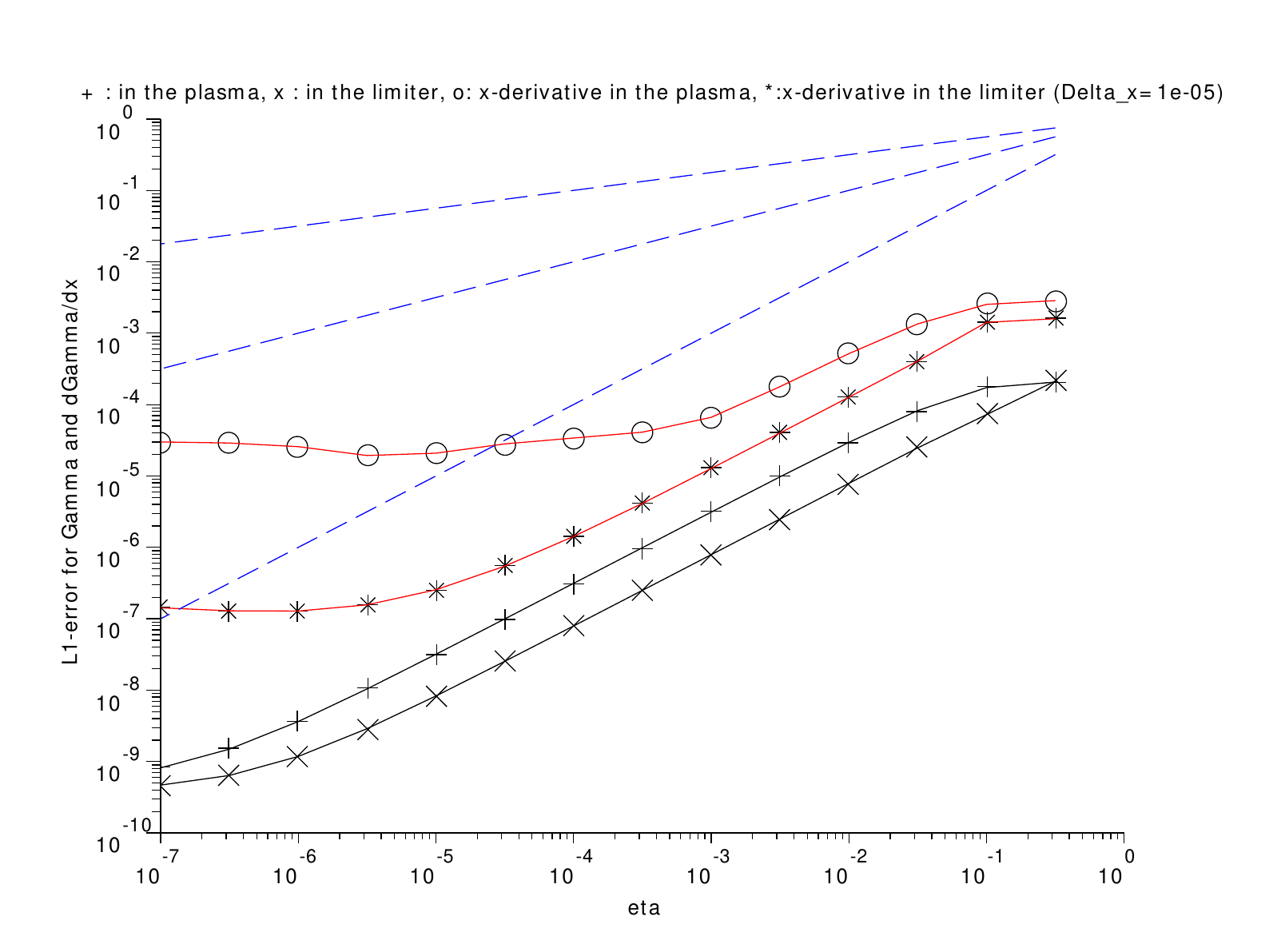}
\end{center}
\caption{
Errors for $N$, $\partial_x N$, $\Gamma$ and $\partial_x \Gamma$ in $L^1$ norms with the penalization free of boundary layer. The dashed lines represent the curves
$\eta^{\frac{1}{4}}, \eta^{\frac12}$ and $\eta$.
}
\label{Error_newpen_regular}
\end{figure}

One of the main difficulties for the implementation of the penalization, is the choice of a boundary condition at $x=0.5$ which is necessary for the numerical scheme. As only $\Gamma$ is penalized, we need a transparent boundary condition for $N$. 
For the numerical tests, the boundary condition comes from the asymptotic expansion up to the first order of the BKW analysis.
We carry out the computations up to $t=1$ with an adaptive time step so that the CFL condition is always satisfied. The results are plotted in Fig. \ref{Plot_newpen_regular}. 
In Fig. \ref{Error_newpen_regular}, we observe that the optimal rate of convergence $\mathcal{O}(\eta)$ is reached for the $L^1$ norms, even for the derivatives. This gives a numerical evidence of the absence of boundary layer. 
The same numerical results in $\mathcal{O}(\eta)$ are obtained if the penalty term in (\ref{Penal_pb_OK}) is replaced by $\frac{\chi}{\eta} \left(\frac{\Gamma}{N}-M_0 \right)$, see \cite{Aup10}.

When the parameter $\epsilon=0.01$, i.e. close to a characteristic boundary, the computations show that, for $\eta$ sufficiently small, $\eta \leq \mathcal{O}(\epsilon)$, the convergence results are similiar; see details in \cite{Ang11}.

\noindent{\bf Acknowledgements:} This work has been funded by the ANR ESPOIR (Edge Simulation of the Physics Of ITER Relevant turbulent transport)and the \emph{F\'ed\'eration nationale de Recherche Fusion par Confinement Magn\'etique} (FR-FCM). We thank Guillaume Chiavassa, Guido Ciraolo and Philippe Ghendrih for fruitful discussions.

\bibliographystyle{Authors/styles/spmpsci.bst}
\bibliography{biblio_article}

\begin{thebibliography}{10}
\providecommand{\url}[1]{{#1}}
\providecommand{\urlprefix}{URL }
\expandafter\ifx\csname urlstyle\endcsname\relax
  \providecommand{\doi}[1]{DOI~\discretionary{}{}{}#1}\else
  \providecommand{\doi}{DOI~\discretionary{}{}{}\begingroup
  \urlstyle{rm}\Url}\fi

\bibitem{Ang11}
Angot, P., Auphan, P., Guès, O.: An optimal penalty method for the hyperbolic
  system modelling the edge plasma transport in a tokamak.
\newblock Preprint in preparation  (2011)

\bibitem{Aup10}
Auphan, T.: Méthodes de pénalisation pour des systèmes hyperboliques
  application au transport de plasma en bord de tokamak.
\newblock Master's thesis, Ecole Centrale Marseille (2010)

\bibitem{ben07}
Benzoni-Gavage, S., Serre, D.: Multidimensional hyperbolic partial differential
  equations. First-order systems and applications.
\newblock Oxford Mathematical Monographs. Oxford University Press (2007)

\bibitem{Bou98}
Bouchut, F., James, F.: One-dimensional transport equations with discontinuous
  coefficients.
\newblock Nonlinear Anal. \textbf{32}, 891--933 (1998)

\bibitem{For08}
Fornet, B.: Small viscosity solution of linear scalar 1-d conservation laws
  with one discontinuity of the coefficient.
\newblock Comptes Rendus Mathematique \textbf{346}(11-12), 681 -- 686 (2008)

\bibitem{For09}
Fornet, B., Guès, .: Penalization approach of semi-linear symmetric hyperbolic
  problems with dissipative boundary conditions.
\newblock Discrete and Continuous Dynamical Systems \textbf{23}(3), 827 -- 845
  (2009)

\bibitem{Gal03}
Gallouët, T., Hérard, J.M., Seguin, N.: Some approximate godunov schemes to
  compute shallow-water equations with topography.
\newblock Computers and Fluids \textbf{32}(4), 479 -- 513 (2003)

\bibitem{gue90}
Guès, O.: Problème mixte hyperbolique quasi-linéaire caractéristique.
\newblock Communications in Partial Differential Equations \textbf{15},
  595--654 (1990)

\bibitem{Iso10}
Isoardi, L., Chiavassa, G., Ciraolo, G., Haldenwang, P., Serre, E., Ghendrih,
  P., Sarazin, Y., Schwander, F., Tamain, P.: Penalization modeling of a
  limiter in the tokamak edge plasma.
\newblock Journal of Computational Physics \textbf{229}(6), 2220 -- 2235 (2010)

\bibitem{Pou97}
Poupaud, F., Rascle, M.: Measure solutions to the linear multi-dimensional
  transport equation with non-smooth coefficients.
\newblock Communications in Partial Differential Equations \textbf{22},
  225--267 (1997)

\bibitem{rau85}
Rauch, J.B.: Symmetric positive systems with boundary characteristic of
  constant multiplicity.
\newblock Trans. Amer. Math. Soc. \textbf{291}(1), 167--187 (1985)

\bibitem{rau74}
Rauch, J.B., Massey, F.J.I.: Differentiability of solutions to hyperbolic
  initial-boundary value problems.
\newblock Trans. Amer. Math. Soc. \textbf{189}, 303--318 (1974)

\bibitem{Tam07}
Tamain, P.: Etude des flux de matière dans le plasma de bord des tokamaks,
  alimentation, transport et turbulence.
\newblock Ph.D. thesis, Université de Provence (2007)

\end{thebibliography}

  \noindent{\small The paper is in final form and no similar paper has been
  or is being submitted elsewhere.}
\end{document}